\documentclass[a4paper,11pt,abstracton=true,headsepline,bibtotoc]{scrartcl}

\RequirePackage{amsmath}                        
\RequirePackage{mathtools}
\RequirePackage{amssymb}                        
\RequirePackage{latexsym}                       
\RequirePackage{stmaryrd}                       

\RequirePackage[all]{xy}                        
\UseTips                                        
\newdir^{c}{{}*!/-3pt/\dir^{(}}
\newdir_{c}{{}*!/-3pt/\dir_{(}}
\newdir^{ (}{{}*!/-3pt/\dir^{(}}
\newdir_{ (}{{}*!/-3pt/\dir_{(}}

\RequirePackage{xspace}                         
\RequirePackage{calc}                           
\RequirePackage{ifthen}                          

\RequirePackage{mathrsfs}
\renewcommand{\mathcal}[1]{\mathscr{#1}}

\RequirePackage[dvipsnames,usenames]{xcolor}  
\RequirePackage{hyperref}
\hypersetup{
bookmarksdepth=3,
bookmarksopen,
bookmarksnumbered,
pdfstartview=FitH,
colorlinks,
linkcolor=RawSienna,
anchorcolor=BurntOrange,
citecolor=OliveGreen,
filecolor=BlueViolet,
menucolor=Yellow,
urlcolor=OliveGreen
}

\usepackage[textsize=footnotesize]{todonotes}
\newcommand{\mb}[1]{\todo[author=MB, color=orange!40]{#1}}
\newcommand{\mbl}[1]{\todo[inline, author=MB, color=orange!40]{#1}}   
\newcommand{\df}[1]{\todo[author=DF, color=green!40]{#1}}

\usepackage{colonequals}
\usepackage[german,english]{babel}      
\usepackage[utf8]{inputenc}
\usepackage{mabliautoref}

\theoremstyle{plain}   
\basetheorem{thm}{Theorem}{}

\maketheorem{prop}{Proposition}{thm}
\maketheorem{cor}{Corollary}{thm}
\maketheorem{lem}{Lemma}{thm}

\theoremstyle{definition}    

\maketheorem{defn}{Definition}{thm}

\theoremstyle{remark}    

\maketheorem{rem}{Remark}{thm}
\maketheorem{ex}{Example}{thm}

\numberwithin{equation}{section}

\title{Cartier Crystals have finite global dimension}
\author{Manuel Blickle and Daniel Fink}
\date{11.11.2022}

\newcommand{\Hom}{\operatorname{Hom}}
\newcommand{\Ext}{\operatorname{Ext}}
\newcommand{\QCoh}{\operatorname{QCoh}}
\newcommand{\uQCoh}{\operatorname{uQCoh}}
\newcommand{\LNil}{\operatorname{LNil}}
\newcommand{\QCrys}{\operatorname{QCrys}}
\newcommand{\Crys}{\operatorname{Crys}}
\newcommand{\colim}{\operatornamewithlimits{colim}}
\newcommand{\id}{\operatorname{id}}
\newcommand{\usc}[1][m]{\underline{\phantom{#1}}}
\newcommand{\Spec}{\operatorname{Spec}}

\renewcommand{\phi}{\varphi}

\renewcommand{\theta}{\vartheta}

\renewcommand{\epsilon}{\varepsilon}

\renewcommand{\to}[1][]{\xrightarrow{\ #1\ }}

\newcommand{\into}[1][]{\xhookrightarrow{\ #1\ }}

\begin{document}

\maketitle
\begin{abstract}
    We show that the category of quasi-coherent Cartier crystals is equivalent to the category of unit Cartier modules on an $F$-finite noetherian ring $R$, and that these equivalent categories have finite global dimension, by showing that every quasi-coherent Cartier crystal has a finite injective resolution. The length of the resolution is uniformly bounded by a bound only depending on $R$. Our result should be viewed as a generalization of the main result of Ma \cite{Ma} showing that the category of unit $R[F]$-modules over a $F$-finite \textbf{regular} ring $R$ has finite global dimension $\dim R + 1$. 
\end{abstract}

\footnotetext{Institut für Mathematik, Johannes Gutenberg Universität Mainz, 55099 Mainz, Germany. Both authors were supported by the Deutsche Forschungsgemeinschaft (DFG) through the Collaborative Research Centre TRR 326 ''Geometry and Arithmetic of Uniformized Structures'' }

\section{Notation and Statement}
Let $p$ be a prime and $R$ an $\mathbb{F}_p$-algebra such that the Frobenius $F : R \to R, r \mapsto r^p$ is a finite map of rings (i.e. $R$ is $F$-finite). As defined in \cite{BBCart}, a \emph{Cartier module} $M$ is a right module over the non-commutative ring 
\[
    R[F] = \frac{R\{F\}}{\langle r^pF-Fr | r \in R \rangle}. 
\]
As a left $R$--module $R[F]\cong \oplus_{e \in \mathbb{N}} RF^e$ is free. We denote the category of Cartier modules (= right $R[F]$-modules) by $\QCoh_{\kappa}(R)$. The right action $\  \cdot F$ on $M$ is encoded in an $R$-linear map $C_M: F_*M \to M$, called the \emph{structural map of the Cartier module}. By adjunction $F_* \dashv F^\flat$, such a map is equivalent to an $R$-linear map $\kappa_M \colon M \to F^\flat M=\Hom_R(F_*R,M)$ where the right hand side is viewed as an $R$ module via the ring isomorphism $R \cong F_*R$ given by the identity. The same holds for powers of $F^e$ yielding the structural maps $C^e_M : F^e_*M \to M$ and $\kappa^e_M : M \to F^{e\flat}M$ with the expected compatibilities.   

A Cartier module $M$ is called \emph{nilpotent} if some sufficiently large power of $F$ acts as zero on $M$. It is called \emph{locally nilpotent} if it is the union of its nilpotent Cartier submodules.   If $R$ is $F$-finite it is shown in \cite[2.1.4]{BBCrys} that the subcategory $\LNil_{\kappa}(R)$ of locally nilpotent Cartier modules is a Serre subcategory of the Grothendieck category $\QCoh_{\kappa}(R)$. Hence the localized category 
\[
\QCrys_{\kappa}(R) := \QCoh_\kappa(R)/\LNil_\kappa(R)
\]
is also a Grothendieck category, called the \emph{category of quasi-coherent Cartier crystals}. In particular $\QCrys_\kappa(R)$ has enough injectives. In \cite{BBCart} it was shown that the subcategory $\Crys_\kappa(R)$ of $\QCrys_\kappa(R)$ consisting of Cartier modules whose underlying $R$-module is finitely generated enjoys remarkable finiteness properties, e.g. it is noetherian and artinian. The point of this note is to show that also homologically Cartier crystals are small:
\begin{thm}\label{thm.main}
For an $F$-finite noetherian $\mathbb{F}_p$-algebra $R$, the category of quasi-coherent Cartier crystals $\QCrys_\kappa(R)$ has finite injective dimension.  
\end{thm}
To explain the connection to Ma's result refered to in the abstract and to set up some notation for our proof below, consider the \emph{unitalization functor}
\[
    u: \QCoh_{\kappa}(R) \to \QCoh_{\kappa}(R), M \mapsto \colim_n F^{n\flat}M 
\]
where the colimit on the right is along the iterates of the adjoint structural maps. Unitalization sends $\LNil_\kappa(R)$ to zero and hence factors through the localization $\QCrys_\kappa(R)$. Its essential image are those Cartier modules $N$ such that the adjoint structural map $\kappa_M : M \to F^\flat M$ is an isomorphism. We call this essential image the category of \emph{unit Cartier modules} denoted $\QCoh_{\kappa}^{\mathrm{unit}}(R)$. 

If $R$ is $F$-finite and noetherian, unitalization is localizing with kernel precisely the locally nilpotent Cartier modules, it induces an equivalence $\QCrys_\kappa(R) \to \QCoh_{\kappa}^{\mathrm{unit}}(R)$, see \autoref{prop.CrysIsUnit} below. In particular we equivalently show, that the category of unit Cartier modules $\QCoh_{\kappa}^{\mathrm{unit}}(R)$ has finite injective resolutions for any $F$-finite noetherian ring $R$.

If $R$ is \textbf{regular} noetherian and $F$-finite, the category of \emph{unit $R[F]$-modules}, which was introduced by Lyubeznik in \cite{Lyub} (see also \cite{EK,BhattLurie}), consists of \textbf{left} $R[F]$-modules (i.e. $R$-modules equipped with a map $F_M: M \to F_*M$), such that the adjoint structural map 
\[
\theta_M: F^* M \to M
\]
is an isomorphism. Tensoring with the canonical module $\omega_R$ (which is invertible as $R$ is regular and $F$-finite) yields an equivalence between the abelian categories of unit $R[F]$-modules and unit Cartier modules, see also \cite{Ma,BBCart}. The main result in \cite[Theorem 1.3]{Ma} asserts that in the case that $R$ is regular noetherian and $F$-finite, this category of unit $R[F]$-modules has finite injective dimension equal to $\dim R + 1$. Hence our \autoref{thm.main} should be considered as a generalization of Ma's result to the case of arbitrary $F$-finite noetherian rings. 

\subsection*{Aknowledgements}
The main result of this paper was part of the master thesis of the second author and builds on unpublished results of the first author with Gebhard Böckle. We thank Gebhard Böckle for graiously allowing to include these in this paper.
\section{The Proof}

The argument consists of 3 steps, each of which uses a certain adjunction (together with some additional information) to be able to control the injectives in the various categories. The first two steps yield, together with the identifications explained above, a self contained and conceptual proof of Ma's result. The third step is our new contribution and is an application of previous still unpublished work of the first author with Böckle in \cite{BBCrys}.

Let $R$ be an $\mathbb{F}_p$-algebra. We will show that:
\begin{enumerate}
    \item If $R$-Mod has finite global dimension $d$ (e.g. $R$ is regular noetherian) then the category of quasi-coherent Cartier Modules $\QCoh_\kappa(R)$ has finite global dimension $d+1$.
    \item For $R$ $F$-finite and regular, the injective dimension of $\QCrys_{\kappa}(R)$ is bounded above by the injective dimension of $\QCoh_\kappa(R)$.
    \item If $R$ is $F$-finite and regular and $A=R/\mathfrak{a}$ is a quotient of $R$, then the injective dimension of  $\QCrys_\kappa(A)$ is bounded above by the injective dimension of $\QCrys_{\kappa}(R)$.
\end{enumerate}
Together with the observation of Gabber \cite[Remark 13.6]{Gab} that every $F$-finite noetherian ring $A$ is the quotient of a regular $F$-finite noetherian ring $R$, these steps show our main result. They also yield that the injective dimension of $A$ is bounded above by $d+1$ where $d=\dim R$ for any regular noetherian $F$-finite $R$ of which $A$ is a quotient. We expect that this bound is optimal.

\subsection*{Proof of part (a)}
First consider the canonical ring map $R \to R[F]$. It yields the adjunction for extension and restriction of scalars
\[
\xymatrix{
\QCoh_\kappa(R) \ar@/^/[rr]^{\operatorname{Rest.}} && R-\operatorname{Mod} \ar@/^/[ll]^{\cdot \otimes_R R[F]}.
}
\]
As a left $R$-module, $R[F] \cong \oplus_{e \in \mathbb{N}} RF^e$ is free, hence the extension of scalars $\cdot \otimes_R R[F]$ is an exact functor from right $R$-modules to Cartier modules. The restriction functor is right adjoint to this exact functor and hence preserves injectives. Let $M \to I^\bullet$ be an injective resolution of the Cartier module $M$. By what was just said, it is also an injective resolution of the underlying $R$-module. Hence the $\Ext_{R[F]}$ groups for Cartier modules $N \otimes_R R[F]$ arising from an $R$ module $N$ by extension
\[
    \Ext^i_{R[F]}(N \otimes_R R[F],M) \cong h^i(\Hom_{R[F]}(N \otimes_R R[F],I^\bullet)) \cong h^i(\Hom_{R}(N,I^\bullet)) = \Ext^i_{R}(N,M)
\]
can be computed as $\Ext_R$ of the $R$-Module $N$.

In \cite[Lemma 3.1]{Ma} a two step resolution of any Cartier Module $N$ 
\[
    0 \to F_*N \otimes_R R[F] \to N \otimes_R R[F] \to N \to 0
\] 
by induced modules is shown, where the right map is the $R[F]$-action and the left is $C_N \otimes 1 - 1\otimes F$. Applying to this short exact sequence the long exact sequence for $\Ext^i_{R[F]}(\cdot,M)$ together with the above identification yields exact sequences
\[
\to \Ext^i_{R}(F_*N,M) \to \Ext^{i+1}_{R[F]}(N,M) \to \Ext^{i+1}_{R}(N,M) \to 
\]
for all $i \geq 0$. Hence if $R$ has finite global dimension $d$, i.e. $\Ext^i_{R}(N,M)=0$ for all $R$-Modules $M,N$ and $i > d$, then $\Ext^i_{R[F]}(N,M)=0$ for all $i>d+1$. Hence the category $\QCoh_{\kappa}(R)$ of quasi-coherent Cartier modules has global dimension $\leq d+1$.

\subsection*{Proof of part (b)}
We employ the unitalization which was mentioned in the introduction. Let $M$ be a Cartier module over $R$ with adjoint structural map $\kappa_M : M \to F^\flat M$. Via $F^\flat (\kappa_M)$ the $R$-module $F^\flat M$ has a natural Cartier module structure such that $\kappa_M$ is a map of Cartier modules. Iterating this we obtain a directed system
\[
    M \to[\kappa] F^\flat M \to[F^\flat\kappa] F^{2\flat}M \to[F^{2\flat}\kappa] F^{3\flat}M \to \ldots .
\]
We define the \emph{unitalization} of $M$ as $u(M) = \colim_e F^{e\flat}(M)$ the colimit (as Cartier modules) of the above system. Since $F^e_*R$ is finitely presented in $R$-Mod if $R$ is Noetherian and $F$-finite, the functor $F^{e\flat}$ commutes with filtered colimits which implies that the adjoint structural map $\kappa_{u(M)} : u(M) \to F^\flat u(M)$ is an isomorphism (its inverse is induced by the identities $F^{e\flat}N \to F^{e\flat}N$ in the limit system). Cartier modules $N$ with this property that the adjoint structural map is an isomorphism are called \emph{unit Cartier modules}. One has the following easy lemma:
\begin{lem}\label{lem.u}
Let $R$ be noetherian and $F$-finite, and $M$ a Cartier module. 
\begin{enumerate}
    \item There is a natural transformation $\eta: \id \to u$. 
    \item $u(\eta) = \eta u: u \to u^2$ is an isomorphism.
    \item $M$ is locally nilpotent if and only if $u(M) = 0$ 
    \item $M$ is a unit Cartier module if and only if $\eta: M \to u(M)$ is an isomorphism.
\end{enumerate}
\end{lem}
\begin{proof} 
The natural transformation of (a) is induced by the natural morphism into the colimit. Assertion (b) holds since $u(\kappa_M)=\kappa_{u(M)}$ is an isomorphism. 

For (c) we observe that if  $u(M)$ vanishes, each element $m\in M$ lies in the kernel of some power of $\kappa_M$, i.e. in the largest Cartier submodule of $M$ on which $F^n$ acts as zero, which is obviously nilpotent. Conversely, since colimits commute with colimits, we can assume that $M$ is nilpotent. In this case, some power of $F$ acts as zero and therefore the same power of $\kappa_M$ is zero. Hence $u(M)$ vanishes since every map $F^{e\flat}M\to u(M)$ is zero. 

(d) If $M$ is a unit Cartier module, each power of $\kappa_M$ is an isomorphism and therefore the same holds for $\eta$. If conversely $\eta$ is an isomorphism, the stuctural map $\kappa_M$ is an isomorphism since all other morphisms in the following diagram are:
\[
\xymatrix{
M \ar[rr]^\eta \ar[d]_{\kappa_M} && u(M) \ar[d]^{\kappa_{u(M)}} \\ 
F^\flat M \ar[rr]_{F^\flat(\eta)} && F^\flat(u(M))
}
\]
\end{proof}
This has the following consequence:
\begin{prop}\label{prop.unitalization}
Let $R$ be noetherian and $F$-finite. The exact localization functor $\pi:\QCoh_\kappa(R) \to \QCrys_\kappa$ has a right adjoint 
\[
\xymatrix{
\QCoh_\kappa(R) \ar@/^/[rr]^{\pi} && \QCrys_\kappa(R) \ar@/^/[ll]^{\overline{u}}
}
\]
for which the unit of adjunction $\pi \circ \overline{u} \to \id$ is an isomorphism (equiv. $\overline{u}$ is fully faithful). If $R$ is further assumed to be regular, then $\overline{u}$ is exact and on underlying Cartier modules given by the unitalization functor $u$ (i.e. $u = \overline{u} \circ \pi$). 
\end{prop}
\begin{proof}
The existence of a fully faithful right adjoint functor is a consequence of  \cite[Proposition 2.4.8]{BPCrys} since $\QCoh_\kappa(R)$ is a Grothendieck category and for $R$ $F$-finite $\LNil_\kappa(R)$ is a the Serre subcategory \cite[2.1.4]{BBCrys} which is closed under filtered colimits.

If $R$ is in addition regular, i.e. the Frobenius of $R$ is flat \cite{KunzFlat}, then $F^{e\flat}$ is an exact functor and hence $u=\colim_e F^{e\flat}$ is exact as a colimit of exact functors. The exactness of $u$ furthermore implies by \autoref{lem.u} (c) that $u$ maps precisely the nil-isomorphism (i.e. morphisms of Cartier modules such that kernel and cokernel are locally nilpotent) to isomorphism so that $u$ factors through $\QCrys_\kappa(R)$. It is clear that the resulting functor is right adjoint to $\pi$, hence equal to $\overline{u}$. It follows from the general theory of localization \cite[Corollary 3.12]{Pop} that since $u$ is exact so is $\overline{u}$.
\end{proof}

\begin{cor}\label{thm.ureflectsInjectivity}
If $I \in \QCrys_\kappa(R)$ is injective then so is $\overline{u}(I) \in \QCoh_\kappa(R)$. If $R$ is regular, then the converse holds as well.  
\end{cor}
\begin{proof}
If $I \in \QCrys_\kappa(R)$ is injective then so is $\overline{u}(I) \in \QCoh_\kappa(R)$ since $\overline{u}$ is right adjoint to an exact functor. The fully faithulness of $\overline{u}$ shows 
\[
    \Hom_{\QCrys}(M,I) \cong \Hom_{\QCoh}(\overline{u}(M),\overline{u}(I)).
\]
If $R$ is regular, $\overline{u}$ is exact and hence $\Hom_{\QCoh}(\overline{u}(M),\overline{u}(I))$ is an exact functor in $M$ provided $\overline{u}(I)$ is injective. Hence if $\overline{u}(I)$ injective, then $\Hom_{\QCrys}(M,I)$ is exact in $M$ and therefore $I$ is injective.  
\end{proof}
\begin{proof}[Proof of (b)]
Let
\[
0 \to M \to I^0 \to \ldots \to I^d \to C \to 0
\]
be an exact sequence in $\QCrys_\kappa(R)$ with $I^j$ injective. As $R$ is regular $F$-finite, $\overline{u}$ is exact and hence we get the exact sequence
\[
0 \to \overline{u}(M) \to \overline{u}(I^0) \to \ldots \to \overline{u}(I^d) \to \overline{u}(C) \to 0.
\]
in $\QCoh_\kappa(R)$ where the $\overline{u}(I^j)$ are injective. If the injective dimension of $\QCoh_\kappa(R)$ is $\leq d$ then it follows by dimension shifting that $\overline{u}(C)$ is injective in $\QCoh_\kappa(R)$. But now it follows that $C$ is injective in $\QCrys_\kappa(R)$ by \autoref{thm.ureflectsInjectivity}. Hence $M$ has an injective resolution of length $\leq d$ which concludes the proof of this step.
\end{proof}
\begin{rem}\label{rem.computingderived}
This proof in particular shows that if $R$ is regular and $F$-finite and $M \in \QCoh_\kappa(R)$ then $u(M)$ has a (finite) injective resolution $uM \to J^\bullet$ in $\QCoh_\kappa(R)$ such that $\pi M  \cong \pi u(M)  \to \pi J^\bullet$ is an injective resolution of $\pi M$ in $\QCrys_\kappa(R)$ (in fact, $J^\bullet = \overline{u}I^\bullet$ for an injective resolution $I^\bullet$ of $M$ in $\QCrys_\kappa(R)$). 

Hence if $G$ is a (left exact) functor between categories of Cartier modules which preserves nil-isomorphism,
and hence descends to a functor $\overline{G}$ on $\QCrys_\kappa$ such that $\overline{G}\circ \pi = \pi \circ G$, the right derived functor $R\overline{G}$ of $\overline{G}$ is given by 
\[
 R\overline{G}(\pi M) = \overline{G}(\pi J^\bullet) = \pi G(J^\bullet)=\pi RG(uM).
\]
where $\pi$ also denotes the derived functor of the exact localization functor $\pi$ from Cartier modules to Cartier crystals.
\end{rem}

Also directly from \autoref{lem.u} it follows that unitalization is left adjoint to the inclusion of the full subcategroy of \emph{unit Cartier modules} $\QCoh_\kappa^{\mathrm{unit}}(R)$ into $\QCoh_\kappa(R)$, the categroy of all Cartier modules. By definition, the essential image of $u$ is $\QCoh_\kappa^{\mathrm{unit}}(R)$, which can be seen to be an Grothendieck abelian category. Unless $R$ is regular, the inclusion $\iota : \QCoh_\kappa^{\mathrm{unit}}(R) \into \QCoh_\kappa(R)$ is only left exact, but not exact (the cokernel in unit Cartier modules is the unitalization of the cokernel in Cartier modules). As a functor $u: \QCoh_\kappa(R) \to \QCoh_\kappa^{\mathrm{unit}}(R)$ unitalization is exact, hence the inclusion $\iota$ preserves injectives. Furthermore, restricting the right hand side of \autoref{prop.unitalization} to $\QCoh_\kappa^{\mathrm{unit}}(R)$ then yields an equivalence of categories:
\begin{prop}\label{prop.CrysIsUnit}
For any $F$-finite noetherian Ring $R$ unitalization induces an equivalence of categories
\[
\xymatrix{
\QCoh_\kappa^\mathrm{unit}(R) \ar@/^/[rr]^{\pi}_\cong && \QCrys_\kappa(R) \ar@/^/[ll]^{\overline{u}}
}
\] 
\end{prop}
\begin{proof}
As $\QCoh_\kappa^{\mathrm{unit}}(R)$ is defined as the full subcategory, it follows from \autoref{prop.unitalization} that the functor $\overline{u}: \QCrys_\kappa(R) \to \QCoh_\kappa^{\mathrm{unit}}(R)$ is still fully faithful and, by \autoref{lem.u} (d), essentially surjective.
\end{proof}

\subsection*{Proof of part (c)}
The proof of part (c), again,  relies on an adjunction, which was obtained by the first author and Böckle as part of a Kashiwara-type equivalence for Cartier Crystals \cite[Section 4.1]{BBCrys}. For completeness we summarize the argument in the special case needed.

Let $R \to A=R/\mathfrak{a}$ be a surjective map of $F$-finite $\mathbb{F}_p$-algebras corresponding to the closed immersion $\Spec R/\mathfrak{a} \into \Spec R$. Denote by $U \subset 
\Spec R$ the open complement of $\Spec R/\mathfrak{a}$ in $\Spec R$. 
\[ 
i: \Spec R/\mathfrak{a} \hookrightarrow \Spec R \hookleftarrow U : j \]
It is verified in \cite[Section 2]{BBCart} that the functors of modules $i_*$, $i^\flat = \Hom_R(R/\mathfrak{a},\usc)$, $j^\flat=j^*$ and $j_*$ naturally extend to functors on Cartier modules. Furthermore, the well known adjunctions $i^* \dashv i^\flat$ and $j^* \dashv j_*$ for modules are indeed adjunctions in the category of Cartier modules. This can be verified by observing that the units of adjunction are Cartier linear, i.e. preserve the right action of Frobenius. For example in the case of $i^* \dashv i^\flat$ this comes down to checking that 
\[
    \eta: i_* i^\flat \to[\operatorname{ev}_1] \id \text{ and } \epsilon: \id \to[m \mapsto (s \mapsto ms)] i^\flat i_*
\]
are morphisms of Cartier modules. But this is obvious as the first one is just the natural inclusion of the $\mathfrak{a}$-torsion submodule $i_*i^\flat M=M[\mathfrak{a}] \subseteq M$ (which is a Cartier submodule) and the second is the identity.

Again by \cite[Section 2.1]{BBCart} each of these functors preserves nilpotence. As $j^*,i_*$ and $j_*$ commute with filtered colimits, these also preserve local nilpotence. As $R$ is noetherian $R/\mathfrak{a}$ is finitely presented over $R$ and therefore $i^\flat = \Hom_R(R/\mathfrak{a},\usc)$ also commutes with filtered colomits. Hence $i^\flat$  preserves local nilpotence as well. One easily concludes \cite[Proposition 3.6]{BBCart} that each of these functors preserve nil-isomorphisms, hence they descend to the respective categories of Cartier crystals. We denote the resulting functors by $\overline{i}_*$, $\overline{i}^*$ and respectively for $j$. 

As the $\mathfrak{a}$-power torsion submodule $\Gamma_\mathfrak{a} M \subseteq M$ is the colimit of Cartier submodules $\colim_n M[\mathfrak{a}^n]$ it is also a Cartier submodule of $M$. Hence $\Gamma_\mathfrak{a}$ is a functor on Cartier modules which preserves local nilpotence. Since injective Cartier modules have injective underlying modules (see part (a)) the right derived functors $Ri^\flat$, $R\Gamma_\mathfrak{a}$ and $Rj_*$ in the category of Cartier modules have as their underlying modules the ``usual'' right derived functors. An important observation is the following lemma.

\begin{lem}\label{lem.loccohom}
The natural map of Cartier modules $i_*i^\flat M \to \Gamma_\mathfrak{a}M$ is a nil-isomorphism. Hence $\pi i_*Ri^\flat M \to \pi R\Gamma_\mathfrak{a}M$ is an isomorphism in $D^+(\QCrys_\kappa(R))$. 
\end{lem}
\begin{proof}
Since the sequence of Frobenius powers $\mathfrak{a}^{[p^e]}$ is cofinal within the usual powers of $\mathfrak{a}$ we have that $\colim_e M[\mathfrak{a}^{[p^e]}] = \Gamma_\mathfrak{a}M$. For $a \in \mathfrak{a}$ and $m \in M$ we have $C_M(a^{p}m)=aC_M(m)$ and hence $C_M(M[\mathfrak{a}^{[p]}]) \subseteq M[\mathfrak{a}]$. This implies that the inclusions $M[\mathfrak{a}^{[p^e]}] \subseteq M[\mathfrak{a}^{[p^{e+1}]}]$ are nil-isomorphisms for all $e$. By taking the colimit this shows the first claim. 

To check the second claim we may represent an object in $D^+(\QCoh_\kappa(R))$ by a quasi-isomorphic bounded below complex of injectives $I^\bullet$. Then $R\Gamma_\mathfrak{a} I^\bullet \cong \Gamma_\mathfrak{a} I^\bullet$ and $i_*Ri^\flat I^\bullet \cong i_*i^\flat I^\bullet$ and the first claim shows that the natural map 
\[
    i_*i^\flat I^\bullet \to \Gamma_\mathfrak{a} I^\bullet
\]
is an object-wise nil-isomorphism, hence yields an isomorphism in $D^+(\QCrys_\kappa(R))$ after applying the exact functor $\pi$. 
\end{proof}
As $R$ is noetherian, one has the standard triangle for local cohomology \cite[IV.§1]{RD}
\begin{equation}\label{eq.triangle} 
    R\Gamma_{\mathfrak{a}} \to \id \to R\Gamma(U,\usc)=Rj_*j^* \to[+1]
\end{equation}
in the bounded below derived category of $R$-modules. By the preceding discussion this is indeed a triangle in $D^+(\QCoh_{\kappa}(R))$ of Cartier modules. Applying $\pi$ we obtain the following result, which is an ad-hoc version of \cite[Theorem 4.1]{BBCrys} in the case that $R$ is regular.
\begin{prop}\label{prop.triangle}
Let $i : R \to R/\mathfrak{a}$ be a surjective map of noetherian $F$-finite rings. Assume (for simplicity) that $R$ is regular. Then there is an exact triangle in $D^+(\QCrys_\kappa(R))$ 
\[
    \overline{i}_*R\overline{i}^\flat \to \id \to \pi Rj_*j^*\overline{u} \to[+1]
\]
where $j$ denotes the open immersion $\Spec R \setminus \Spec R/\mathfrak{a} \into \Spec R$ and $\overline{u}$ is the exact unitalization functor from \autoref{prop.unitalization}. 
\end{prop}
\begin{proof}
Let $M$ be a Cartier Module on $R$. Applying the exact localization functor $\pi$ to the local cohomology triangle above for the Cartier module $uM$ and using the nil-isomorphism of \autoref{lem.loccohom} we obtain the triangle
\[
    \pi i_* Ri^\flat uM \to \pi uM \to \pi Rj_*j^*uM \to[+1]
\]
in $D^+(\QCrys_\kappa(R))$. By definition of $\overline{i}_*$ we have $\pi i_* \cong \overline{i}_* \pi$. By \autoref{rem.computingderived} we have $\pi Ri^\flat u(M) \cong R\overline{i}^\flat \pi M$. As $M \to uM$ is a nil-isomorphism, $\pi M \to \pi uM$ is an isomorphism in $\QCrys\kappa (R)$, the middle term is isomorphic to $\pi M$ in $\QCrys_{\kappa}(R)$. By definition, $\overline{u} \pi M = uM$ which shows the claimed triangle for $\pi M$. But this is sufficient since $\pi$ is essentially surjective.  
\end{proof}
As a corollary we obtain:
\begin{prop}\label{prop.adjcrys}
Let $i : R \to R/\mathfrak{a}$ be a surjective map of noetherian $F$-finite rings. Assume (for simplicity) that $R$ is regular. Then there is an adjunction 
\[
\xymatrix{
  \QCrys_\kappa(R/\mathfrak{a}) \ar@/^/[rr]^{\overline{i}_*} && \QCrys_\kappa(R) \ar@/^/[ll]^{\overline{i}^\flat }
}
\]
and the derived unit of adjunction $\epsilon :N   \to R\overline{i}^\flat \overline{i}_* N$ is a quasi-isomorphism.
\end{prop}
\begin{proof}
The adjunction on the level of Cartier modules was discussed above at the beginning of step (c). As $i_*$ and $i^\flat$ preserve nil-isomorphisms, the adjunction passes to an adjunction of $\overline{i}_*$ and $\overline{i}^\flat$ for Cartier crystals. 

Let $N \in \QCoh_\kappa(R/\mathfrak{a})$ be a Cartier module. The triangle in \autoref{prop.triangle} for $\overline{i}_*\pi N$ has as its third term $Rj_*j^*u i_*N$. For an open immersion $j^*=j^\flat$ and hence $j^*$ commutes with $F^\flat$. Therefore $j^*$ commutes with $u = \colim_e F^{e\flat}$. Since already on underlying quasi-coherent $\mathcal{O}$-modules $j^*i_* = 0$ it follows that the last term is zero. Hence the triangle degenerates to a quasi-isomorphism
\[
    \eta \overline{i}_* : \overline{i}_*R\overline{i}^\flat \overline{i}_*\pi N \to[\cong] \overline{i}_*\pi N .
\]
The adjunction formalism implies that this is the inverse of $\overline{i}_*\epsilon$. Since $\overline{i}_*$ reflects quasi-isomorphisms it follows that $\epsilon$ itself is a quasi-isomorphism.
\end{proof}

\begin{proof}[Proof of (c)]
Let $N$ be in $\QCrys_\kappa(R/\mathfrak{a})$ and let 
\[
    \overline{i}_*N \cong I^\bullet
\]
be a finite injective resolution of $\overline{i}_*N$ in $\QCrys_\kappa(R)$ which exists by step (b) since $R$ is regular. By definition of the right derived functor we have $R\overline{i}^\flat \overline{i}_*N \cong \overline{i}^\flat I^\bullet$. Composing this with the (the inverse of the) quasi-isomorphism $N \cong R\overline{i}^\flat \overline{i}_*N$ of \autoref{prop.adjcrys} we obtain a resolution of $N$
\[
    N \cong \overline{i}^\flat I^\bullet
\]
As the right adjoint of an exact functor $\overline{i}^\flat$ preserves injectives, so that $\overline{i}^\flat I^\bullet$ is a bounded complex of injectives in $\QCrys_\kappa(R/\mathfrak{a})$ of length at most that of $I^\bullet$. This concludes our proof.
\end{proof}

\subsection{Concluding remark}
Rather formally one may derive the following statement for the derived categories of Cartier crystals which could be viewed as weak analog of \cite[12.4.1]{BhattLurie} for Frobenius modules. 
\begin{thm}
Let $R$ be an $F$-finite noetherian $\mathbb{F}_p$-algebra. For any $M^\bullet \in D^b(\QCoh_\kappa(R))$ there is a morphism of complexes of Cartier modules
\[
    M^\bullet \to I^\bullet
\]
such that
\begin{enumerate}
    \item $I^\bullet$ is a bounded complex of injectives in $\QCoh_\kappa(R)$, in particular the $I^j$'s are injective $R$-modules.
    \item $\pi M^\bullet \to \pi I^\bullet$ is a quasi-isomorphism in $D^b(\QCrys_\kappa(R))$
\end{enumerate}
\end{thm}
\begin{proof}
Let $M^\bullet \to u(M^\bullet)$ be the natural map of $M^\bullet$ to its unitalization. $u(M^\bullet)$ is a bounded complex of unit Cartier modules. As all unit Cartier modules have finite injective resolutions by \autoref{prop.CrysIsUnit} and \autoref{thm.main} it follows that there is a map of complexes
\[
    u(M^\bullet) \to I^\bullet
\]
which is a quasi-isomorphism and where $I^\bullet$ is a bounded complex of injectives in the category $\QCoh_\kappa^\mathrm{unit}(R)$ (by hand one can construct $I^\bullet$ splitting $u(M^\bullet)$ into short exact sequences and, using the horseshoe lemma, construct a double complex of injectives whose columns are resolutions of the $u(M^j)$'s in $\QCoh_\kappa^\mathrm{unit}(R)$. The total complex of this double complex is then $I^\bullet$). Since the inclusion $\QCoh_\kappa^\mathrm{unit}(R) \into \QCoh_\kappa(R)$ preserves injectives, the $I^j$'s are indeed injective in $\QCoh_\kappa(R)$. This shows that the composition
\[
    M^\bullet \to u(M^\bullet) \to[\simeq] I^\bullet. 
\]
satisfies (a). Applying the exact localization functor $\pi$ yields (by \autoref{prop.CrysIsUnit}) (b) since $\pi M^\bullet \to \pi u(M^\bullet)$ is (even object-wise) an isomorphism in $\QCrys_\kappa(R)$.  
\end{proof}



\bibliographystyle{alpha}
\bibliography{references}

\end{document}